\documentclass[9pt]{amsart}

\usepackage{moreverb}

\usepackage[dvips,colorlinks,bookmarksopen,bookmarksnumbered,citecolor=red,urlcolor=red]{hyperref}
\usepackage{color}
\usepackage[T1]{fontenc}

\newtheorem{Thm}{Theorem}
\newtheorem{Prop}[Thm]{Proposition}

\newtheorem{Lem}{Lemma}

\def\Remark{\medskip\noindent {\bf Remark.}\ \ignorespaces}
\def\xd {\mathrm{d}}

\def\e{\varepsilon}
\def\R{{\mathbb R}}
\def\N{{\mathbb N}}

\def\S{{\mathbb S}}

\def\fini{f_i^{\mbox{\tiny in}}}

\def\cini{c_i^{\mbox{\tiny in}}}

\DeclareMathAlphabet{\mathitbf}{OML}{cmm}{b}{it}
\newcommand{\na}{\nabla}
\newcommand{\pa}{\partial}

\newcommand{\dij}{\textrm{\DH}_{ij}}
\newcommand{\dji}{\textrm{\DH}_{ji}}

\newcommand\BibTeX{{\rmfamily B\kern-.05em \textsc{i\kern-.025em b}\kern-.08em
T\kern-.1667em\lower.7ex\hbox{E}\kern-.125emX}}

\begin{document}


\title{On the Maxwell-Stefan diffusion limit for a mixture of monatomic gases}

\author{Harsha Hutridurga and Francesco Salvarani}

\begin{abstract}
Multi-species Boltzmann equations for gaseous mixtures, with analytic cross sections and under Grad's angular cutoff assumption, are considered under diffusive scaling.
In the limit, we {formally} obtain an explicit expression for the binary diffusion coefficients in the Maxwell-Stefan equations.
\end{abstract}


\maketitle


\vspace{-6pt}

\section{Introduction} \label{section_intro}

The macroscopic description of diffusive phenomena in mixtures goes back to the 19th Century, thanks to the work of Maxwell \cite{max1866} and Stefan
\cite{ste1871}, who introduced a coupled system of cross-diffusion equations, nowadays known as the Maxwell-Stefan system.
Since then, a wide literature about multicomponent diffusive phenomena has been published by the Chemical Engineering community (we refer to \cite{kri-wes-97} for a review on the state-of-the-art on the subject), since multicomponent diffusion in fluids plays a crucial role in many (bio)chemical processes. 

However, the mathematical and numerical study of the Maxwell-Stefan system is relatively new, and solid results on the subject appeared only very recently (see, for example, \cite{gio, ern-gio, ern-gio2, gio_book, bot-11, bou-goe-gre, jun-ste-13, mcl-bou-14, bou-gre-sal-12, che-jun-15, bou-gre-pav-16}).

In particular, \cite{bou-gre-sal-15} is devoted to the formal derivation, under the standard diffusive scaling, of the Maxwell-Stefan
diffusion equations from the non-reactive elastic Boltzmann system for monatomic gaseous mixtures, in the vanishing Mach and Knudsen numbers limit.
The approach follows the research line introduced in \cite{bar-gol-lev-89, bar-gol-lev-91, bar-gol-lev-93} and allows to deduce asymptotically a compressible dynamics and to obtain, in the limit, an explicit expression of the binary diffusion coefficients in the Maxwell-Stefan system, depending on the reduced mass of the species, on the temperature and on the cross sections of the kinetic model. 
A peculiar feature of {this} approach is the possibility of obtaining explicit coefficients which could be compared with experimental results and contribute to a better understanding of the quantitative behaviour of gaseous mixtures.

The main assumptions of  \cite{bou-gre-sal-15} are the following:
\vskip2mm
\begin{itemize}
\item the initial conditions are well prepared, and are given by local Maxwellians, all with the same temperature;
\item the process is supposed to be isothermal;
\item the cross sections of the Boltzmann system are of Maxwellian type.
\end{itemize}
\vskip2mm

The first two assumptions are very natural, since the presence of temperature gradients may introduce transport phenomena which could hinder the diffusion process, as shown in \cite{bou-gre-pav-sal}. On the other hand, the third one is not completely satisfactory.
Indeed, even if the hypothesis of Maxwellian cross-sections has been very popular as it could lead to many explicit calculations (Maxwell and Boltzmann themselves used them very often \cite{vil-book}), it has been subsequently noticed that the influence of the collision kernel (or equivalently, of the cross-section) on the solutions to the Boltzmann equation is far from being negligible \cite{vil-book}.

For this reason, it is important to recover the Maxwell-Stefan diffusion equations from the non-reactive elastic Boltzmann equations for mixtures with more realistic cross sections. In this article, we introduce a new strategy which allows us to handle {at the formal level} general factorized cross sections, in which the kinetic collision kernel is an analytic function of its argument and the angular collision kernel is even and satisfies Grad's cutoff assumption \cite{gra-book}.
By considering these more physical cross sections, we are able to obtain explicit quantitative expressions of the binary diffusion coefficients, whose structure exhibits some nonlinearities which were absent in the case of Maxwellian cross sections.

This study is also interesting from a physical point of view. Indeed, even though many results are available in the case of binary mixtures (see, for example, \cite{van1968determination}), many attempts have been made to experimentally determine the binary diffusion coefficients in the Maxwell-Stefan equations, but with a limited success \cite{leahy2005measurements}.
We hope that our computation will give a contribution for encouraging the experimental activity and be useful for comparison with experimental data.

We finally point out that the results obtained in this article are valid for monatomic ideal gases only. The extension of our computations to polyatomic gases should be possible by starting from kinetic systems taking into account the internal energy of the molecules (such as the model proposed in \cite{bourgat1994microreversible}). We will consider this extension in a near future.

The structure of the article is the following: after describing the kinetic model for non-reactive monatomic gaseous mixtures in Section \ref{s2} and 
the Maxwell-Stefan system in Section \ref{section_maxste}, in Section \ref{sn:ms} we compute the diffusion limit and obtain {an} explicit expression of the binary diffusion coefficients in the case of analytical cross sections under Grad's cutoff assumption. Finally, we give an appendix where we explicitly compute the Gaussian integrals in the expression for the binary diffusion coefficients.

\section{The kinetic model}
\label{s2}

The mathematical form of the Boltzmann system for non-reactive monatomic gas mixtures is classical, and has been the starting point of many extensions (such as, for example, \cite{bourgat1994microreversible, des-mon-sal}). However, in order to make this article self-consistent, we briefly describe its form.
 
The model considers a mixture of ideal monatomic inert gases $\mathcal{A}_i$, $i=1,\dots, I$ with $I\ge2$. Each of them is described by
a distribution function $f_i$, which depends on $t\in \R^*_+$ (time),
$x\in\R^3$ (space position) and $v\in \R^3$ (velocity).

By supposing that no chemical reactions occur and no external forces act on the mixture, its time evolution is a consequence of
the mechanical collisions between molecules, which are supposed here to be elastic.

Let us consider two particles belonging to the species $\mathcal A_i$ and $\mathcal A_j$, $1\le i,j\le I$, with masses
$m_i$, $m_j$, and pre-collisional velocities $v'$, $v_*'$. A microscopic collision is an instantaneous phenomenon which modifies the velocities of the particles, which
become $v$ and $v_*$, obtained by imposing the conservation of both momentum and kinetic
energy:
\begin{equation} \label{momen}
m_i v'+m_j v'_* = m_i v+m_j v_*,\qquad 
\frac 12 m_i\,|v'|^2 + \frac 12 m_j\,|v_*'|^2 = \frac 12 m_i\,|v|^2 + \frac 12
m_j\,|v_*|^2.
\end{equation}
The previous equations allow to write $v'$ and $v'_*$ with respect to $v$
and $v_*$:
\begin{equation} \label{v*}
v'=\frac{1}{m_i +m_j}(m_i v+m_j v_* + m_j \vert v- v_*\vert\, \sigma), \qquad 
v'_*= \frac{1}{m_i+m_j}(m_i v+m_j v_* - m_i \vert v- v_*\vert\, \sigma),
\end{equation}
where $\sigma\in\S^2$ describes the two degrees of freedom in (\ref{momen}).

If $f$ and $g$ are nonnegative functions, the operator describing the
collisions of molecules of species $\mathcal A_i$ with molecules of species
$\mathcal A_j$ is defined by
\begin{equation} \label{Q_el_bi}
Q_{ij}(f,g)(v):= \int_{\R^3} \int_{\S^2} B_{ij}(v,v_*,\sigma)\big[ f(v')g(v'_*) -
f(v)g(v_*) \big]  \,\xd\sigma\, \xd v_*,
\end{equation}
where $v'$ and $v'_*,$ are given by \eqref{v*}, and the cross section
$B_{ij}$ satisfies the microreversibility assumptions
$B_{ij}(v,v_*,\sigma)=B_{ji}(v_*,v, \sigma)$ and
$B_{ij}(v,v_*,\sigma)=B_{ij}(v',v'_*,\sigma)$. 

It is clear that, when $i=j$, the previous expressions reduce to the standard Boltzmann kernel in the
mono-species case:
\begin{equation} \label{Q_el_mono}
Q_{ii}(f,f)(v) = \int_{\R^3} \int_{\S^2} B_{ii}(v,v_*,\sigma)\big[ f(v')f(v'_*) -
f(v)f(v_*) \big]  \,\xd\sigma\, \xd v_*.
\end{equation}

The operators $Q_{ij}$ can be written in weak form. For example, by using the changes of variables $(v,v_*) \mapsto (v_*,v)$ and $(v,v_*)
\mapsto (v',v'_*)$, we have
\begin{multline} \label{weak_bi_utile}
\int_{\R^3}Q_{ij}(f,g)(v)\,\psi(v)\, \xd v\\
= -\frac 12 \int_{\R^6} \int_{\S^2} B_{ij}(v,v_*, \sigma)\big[ f(v')g(v'_*) -
f(v)g(v_*) \big] \big[ \psi(v') -\psi(v) \big] \, \xd\sigma\,\xd v\, \xd v_*\\ 
=\int_{\R^6} \int_{\S^2} B_{ij}(v,v_*, \sigma)\,f(v)g(v_*)\, 
\left[ \psi(v') -\psi(v) \right] \, \xd\sigma\,\xd v\, \xd v_*,
\end{multline}
or
\begin{multline} \label{weak_el_bi} 
\int_{\R^3}Q_{ij}(f,g)(v)\,\psi(v)\, \xd v + \int_{\R^3}Q_{ji}(g,f)(v)\, \phi(v)\, \xd v=\\
-\frac 12 \int_{\R^6}\int_{\S^2}B_{ij}(v,v_*,\sigma) \big[ f(v')g(v'_*) -
f(v)g(v_*) \big ] \big[\psi(v') +
\phi(v'_*)-\psi(v)-\phi(v_*)\big]\,\xd\sigma\,\xd v\, \xd v_*,
\end{multline}
for any $\psi$, $\phi: \, \R^3\to \R$ such that the first integrals in
(\ref{weak_bi_utile})--(\ref{weak_el_bi}) are well defined. 

When $\psi(v) = 1$ in (\ref{weak_bi_utile}) we deduce the
conservation of the total number of molecules of species $\mathcal A_i$.
Moreover, if $\psi(v) = m_i\,v$ and $\phi(v_*) = m_j\,v_*$, and then if
$\psi(v) =m_i\,|v|^2/2$ and $\phi(v)=m_j\,|v_*|^2/2$, we recover the conservation of the total momentum and of the total kinetic
energy during the collision between a particle of species $\mathcal A_i$ and a particle of species $\mathcal A_j$:
\begin{equation}\label{consma}
\int_{\R^3} Q_{ij}(f,g)(v) \, \left(
\begin{array}{c}
m_i\,v\\ m_i\,{|v|^2}/2
\end{array}
\right) \, \xd v
+ \, \int_{\R^3}Q_{ji}(g,f)(v) \, \left(
\begin{array}{c}
m_j\,v \\ m_j\,{|v|^2}/2
\end{array}
\right) \, \xd v= 0.
\end{equation}
The system of equations satisfied by the set of distribution functions
$(f_i)_{1\le i\le I}$ is hence
\begin{equation} \label{eqbo_nonscaled}
\partial_t f_i+v \cdot \nabla_x f_i
=  \sum_{j=1}^I Q_{ij}(f_i,f_j) 
\quad \mbox{on } \R_+\times\R^3\times\R^3.
\end{equation}

\section{The Maxwell-Stefan model} \label{section_maxste}

The Maxwell-Stefan model is suitable to describe an ideal gaseous mixture of $I\ge 2$ species, with molecular
masses $m_i$, in which no convective phenomena take place and the system is driven to equilibrium by pure diffusion.

For each species of the mixture $\mathcal{A}_i$,  $1\le i\le I$, we consider its
concentration $c_i$ and its flux $F_i$ -- which depend on the macroscopic variables
$t\in\R^+$ (time) and $x\in\R^3$ (position).

These quantities, which are the unknowns of the system, satisfy the continuity equation
\begin{equation} \label{e:conci}
\pa_t c_i + \na_x\cdot F_i = 0 \quad \mbox{on } \R_+\times \R^3,
\end{equation}
for all $1\le i\le I$.

Let $c=\sum c_i$ be the total concentration of the mixture and let $n_i=c_i/c$
the mole
fraction of species $\mathcal{A}_i$. The Maxwell-Stefan equations can be written in the following form:
\begin{equation} \label{e:cmaxstei}
-c\, \na_x n_i = \frac1c\sum_{j\neq i} \frac{c_j F_i-c_i F_j}{\dij} \quad
\textrm{ on } \R_+\times \R^3,
\end{equation}
for all $1\le i\le I$.

The quantities $\dij$ are the binary diffusion coefficients between the
species $\mathcal{A}_i$ and $\mathcal{A}_j$. They are symmetric with respect to the particles
exchange, in such a way that $\dij=\dji$. 

Note that the Maxwell-Stefan equations \eqref{e:cmaxstei} are linearly dependent. Indeed, by summing \eqref{e:cmaxstei} with respect to $i$, we obtain an identity. Hence, a supplementary equation is necessary for assuring the closure of the Maxwell-Sterfan system \eqref{e:conci}--\eqref{e:cmaxstei}. 

By assuming that the system is closed and under constant and uniform temperature and pressure,
it is usual to assume that the total diffusive flux satisfies 
\begin{equation} \label{e:sumflux}
\sum_{i=1}^I F_i = 0 \quad \mbox{on } \R_+^*\times \R^3,
\end{equation}
which physically means that the diffusive fluxes do not create any mass \cite{gio_book}.

By summing \eqref{e:conci} with respect to $i$, we note that $c$ is uniform in time. Hence, if we suppose that the
molecules of the mixture are initially uniformly distributed, the quantity $c$ is a pure constant.

\section{The Maxwell-Stefan asymptotics} \label{sn:ms}

In this section we apply the strategy proposed in \cite{bou-gre-sal-15}, with {more general} cross sections, satisfying the hypotheses defined in Subsection \ref{cc}.

We work in the standard diffusive scaling, by supposing that the mean free path (or, equivalently, the Knudsen number) tends to zero.
Since temperature gradients can induce transport phenomena, we will suppose that the temperature $T$ of the mixture is a constant.
Moreover, we assume that the bulk velocity of the mixture is of the same order of magnitude as the Knudsen and Mach numbers.

\subsection{Collision kernels}
\label{cc}

In this article, we suppose that the collision kernels $B_{ij}$ depend only on the modulus of the relative velocity, i.e. $|v-v_*|$ and on the cosine of the deviation angle $\theta$, where
$$
\cos\theta=\frac{v-v_*}{|v-v_*|}\cdot\sigma.
$$
{More specifically}, we work with collision kernels of the form:
\begin{align}\label{eq:B-ij}
B_{ij}(v,v_*,\sigma) = \Phi(|v-v_*|) b_{ij}(\cos\theta),
\end{align}
where we assume that the angular collision kernels $b_{ij}\in L^1(-1,+1)$ and are even. Observe that, because of the microreversibility assumption on the collision kernels, we have
\begin{align*}
B_{ij}(v,v_*,\sigma)=B_{ji}(v_*,v,\sigma)\implies b_{ij}\left(\frac{v-v_*}{|v-v_*|}\cdot\sigma\right) = b_{ji}\left(\frac{v_*-v}{|v-v_*|}\cdot\sigma\right),
\end{align*}
and that, by parity, $b_{ij}(\cos\theta) = b_{ji}(\cos\theta)$.

For $\ell\in\{1, 2, 3\}$, we denote by $w_{(\ell)}$ the $\ell$-th component of any vector $w\in\R^3$.

If we introduce the polar variable $\varphi\in [0,2\pi]$, we can find the relationships between the
Euclidean coordinates of $\sigma$ and the spherical ones, namely 
$$
\sigma_{(1)}=\sin \theta\cos\varphi, \quad \sigma_{(2)}=\sin
\theta\sin\varphi, \quad \sigma_{(3)}=\cos \theta.
$$
The kinetic collision kernel $\Phi(|v-v_*|)$ is assumed to be analytic in the following sense: there exists a family $\{a_n\}_{n\in\mathbb{N}^*}\subset\R$ such that $\Phi$ can be written as a uniformly converging even power series:
\begin{align}\label{eq:Phi-power}
\Phi(|v-v_*|) = \sum_{n\in\mathbb{N}^*} a_n |v-v_*|^{2n}.
\end{align}

\subsection{Scaled equation}

In order to arrive at the diffusive limit, we introduce a scaling parameter $0<\e\ll1$ which represents the mean free path.
We denote the corresponding unknown distribution functions as $(f_i^{\e})_{1\le i\le I}$. Each distribution function $f_i^\e$ solves the following scaled version of \eqref{eqbo_nonscaled}:
\begin{equation} \label{eqbo}
\e\, \partial_t f_i^{\e}+v \cdot \nabla_x f_i^{\e}
=  \frac 1\e
 \sum_{j=1}^I Q_{ij}(f_i^{\e},f_j^{\e}), 
\quad \mbox{on } \R_+\times\R^3\times\R^3.
\end{equation}
Finally, we define the corresponding concentrations $(c_i^\e)_{1\le i\le I}$ as the zero-th order moment of the distribution functions $f_i^\e(t,x,v)$:
\begin{equation*}
c_i^\e(t,x)=\int_{\R^3} f_i^\e(t,x, v) \, \xd v, \quad \mbox{ for }(t,x)\in\R_+\times\R^3;
\end{equation*}
this relationship is a fundamental link between the kinetic equations and the Maxwell-Stefan description.

\subsection{Ansatz}

{As we are interested in pure diffusion dynamics, we suppose that the initial data $f^{in}(x,v)$ for the multi-species Boltzmann equations (\ref{eqbo}) are such that
\begin{align*} 
\int_{\R^3}\fini(x,v)\, {\rm d}v = \cini(x),
\qquad
\int_{\R^3} v\fini(x,v)\, {\rm d}v = \mathcal{O}(\e),
\end{align*}
}where
$$
\cini\, :\,\R^3\to \R_+
$$
are $\e$-independent. We moreover suppose that
$$
\sum_{i=1}^I \cini=1 \quad\mbox{on }\R^3, 
$$
which of course implies that each $\cini$ lies in $[0,1]$. {As in \cite{bou-gre-sal-15}, we assume that the evolution following \eqref{eqbo} keeps the distribution functions $f^\e_i(t,x,v)$ in the local Maxwellian state, with a homogeneous temperature $T$}. We hence suppose that there exist 
$$
c_i^{\e}\, :\,\R_+\times \R^3\to \R_+,\qquad
u_i^{\e}\, :\,\R_+\times \R^3\to \R^3,\qquad
1\le i\le I,
$$
such that
\begin{equation} \label{ansatz} 
f^{\e}_i(t,x,v) = {c^{\e}_i(t,x)} \, \left (\frac{m_i}{2\pi k\, T}\right )^{3/2}
e^{-{m_i}\vert v-\e u^{\e}_i(t,x) \vert^2/2kT}, \quad \mbox{ for }(t,x,v)\in\R_+\times\R^3\times\R^3.
\end{equation}
The starting point of our analysis is the following result, proved in \cite{bou-gre-sal-15}:
\begin{Prop}\label{prop:M1-M2}
Under the assumption \eqref{ansatz} on the distribution functions $f^\e_i(t,x,v)$, we have the following mass balance equations for all $1\le i\le I$:
\begin{equation}
\label{M1}
\partial_tc_i^{\e} +\nabla_x\cdot ( c_i^{\e}u_i^{\e})=0\qquad \mbox{ on }\R_+\times\R^3.
\end{equation}
We further have the following momentum balance for all $1\le i\le I$:
\begin{equation}\label{M2}
\e^2 
\left[\pa_t \left( c_i^{\e}u_i^{\e}\right)
+
\nabla_x\cdot\left(c_i^{\e}u_i^{\e}\otimes u_i^{\e}\right)\right] + 
\frac {kT}{m_i}\nabla_x c_i^{\e} = \Theta^\e_i \quad \mbox{ on }\R_+\times\R^3,
\end{equation}
where the $\ell$-th component of $\Theta^\e_i$ is given by
\begin{equation}\label{eq:prop_Theta}
\left(\Theta^\e_i\right)_{(\ell)} =
\frac1\e \sum_{j\neq i} \frac{m_j}{m_i+m_j}\int_{\R^{6}}
\int_{\S^2}
B_{ij}(v, v_*, \sigma)f_i^{\e}(v)f_j^{\e}(v_*)\,
\left({v_*}_{(\ell)}-v_{(\ell)}+ |v-v_*|\sigma_{(\ell)}\right)\, \xd\sigma\,\xd
v_*\, \xd v.
\end{equation}
\end{Prop}

{
\Remark {\it In Equation \eqref{M2}, we have stressed the actual order in $\e$ of various terms. In particular, as shown in the next lemma,
the quantities $\Theta^\e_i$ are of order $\mathcal{O}(1)$. This is of the same order as $\nabla_x c^\e_i$ in \eqref{M2}. 
On the other hand, the first two terms on the left hand side of \eqref{M2} are of order $\mathcal{O}(\e^2)$. 
This clarifies how the diffusive scaling acts on various macroscopic quantities associated with the solutions to the Boltzmann system \eqref{eqbo}.}
}
\bigskip

Our next task is to analyze the right hand side $\Theta^\e_i$ of the momentum balance \eqref{M2}. Observe that $\Theta^\e_i$ depends on the independent variables $(t,x)$. 
\begin{Lem}\label{lem:Theta}
The $\Theta^\e_i(t,x)$ term in the momentum balance \eqref{M2} can be asymptotically approximated as follows:
\begin{align}\label{eq:lem-Theta}
\left(\Theta^\e_i\right)_{(\ell)}(t,x) = \sum_{j\not=i} \Delta_{ij}\left(c^{\e}_i c^{\e}_j (u_j)^{\e}_{(\ell)} - c^{\e}_i c^{\e}_j (u_i^{\e})_{(\ell)}\right) + \mathcal{O}(\e),
\end{align}
where $\Delta_{ij}$ are given by
\begin{equation}\label{eq:Delta-ij}
\begin{array}{cl}
\displaystyle\Delta_{ij} = &\displaystyle a_0 \frac{2\pi  m_j  \|b_{ij}\|_{L^1}}{(m_i+m_j)} +  \, a_1  \frac{10\pi kT \|b_{ij}\|_{L^1}}{m_i}+ \sum_{n\ge2} a_n \frac{ m_j (m_i m_j)^{3/2} \|b_{ij}\|_{L^1}}{4\pi^2 (m_i+m_j)(k\, T)^3}
\int_{\R^{6}} \Big\{\Big( \sum_{n_1+n_2+n_3=n}\frac{n!}{n_1! n_2! n_3!}\times \\[0.2 cm]
&\displaystyle  \prod_{1\le r\le 3}\Big(\sum_{\alpha+\beta=2n_r} \frac{(2n_r)!}{\alpha! \beta!} \left([{v}_{(r)}]^\alpha [-v_{*(r)}]^\beta\right)\Big)\Big)
e^{-{m_i}\vert v \vert^2/2kT}
e^{-{m_j}\vert v_* \vert^2/2kT}\Big\}
\, \xd v_*\, \xd v.
\end{array}
\end{equation}

\end{Lem}

\noindent{\it Proof}\par

\noindent For readers' convenience, let us rewrite the expression \eqref{eq:prop_Theta} for $\left(\Theta^\e_i\right)_{(\ell)}$:
\begin{align*}
\left(\Theta^\e_i\right)_{(\ell)} =
\frac1\e \sum_{j\neq i} \frac{m_j}{m_i+m_j}\int_{\R^{6}}
\int_{\S^2}
B_{ij}(v, v_*, \sigma)f_i^{\e}(v)f_j^{\e}(v_*)\,
\left({v_*}_{(\ell)}-v_{(\ell)}+ |v-v_*|\sigma_{(\ell)}\right)\, \xd\sigma\,\xd
v_*\, \xd v.
\end{align*}
The term containing $\sigma_{(\ell)}$ in \eqref{eq:prop_Theta} vanishes.
Indeed, both terms for $\ell=1$ or $2$ are zero because
$$\int_0^{2\pi}\sin\varphi\,\xd\varphi=\int_0^{2\pi}\cos\varphi\,\xd\varphi=0,$$
and for $\ell=3$, because $b_{ij}$ is even, one has
$$\int_{\S^2} b_{ij}\left(\frac{v-v_*}{|v-v_*|}\cdot\sigma\right)\sigma_{(3)} \,
\xd\sigma=2\pi\int_0^\pi\sin\theta\cos\theta\,
b_{ij}(\cos\theta)\,\xd\theta=2\pi\int_{-1}^1 \eta \,b_{ij}(\eta)\,\xd\eta=0.
$$
Gathering the remaining part of the expression for $\left(\Theta^\e_i\right)_{(\ell)}$, we have:
\begin{align}\label{eq:theta-1}
\left(\Theta^\e_i\right)_{(\ell)} =
\frac{1}{\e} \sum_{j\neq i} \frac{2\pi m_j \|b_{ij}\|_{L^1}}{m_i+m_j}
\int_{\R^{6}}
\Phi(|v-v_*|) f_i^{\e}(v)f_j^{\e}(v_*)\,
\left({v_*}_{(\ell)}-v_{(\ell)}\right)\, \xd
v_*\, \xd v.
\end{align}
We now substitute the ansatz (\ref{ansatz}) for $f_i^\e(v)$ and $f_j^\e(v_*)$ in \eqref{eq:theta-1}.
By writing the power series expansion \eqref{eq:Phi-power} for the kinetic collision kernel we obtain:
\begin{align*}
\left(\Theta^\e_i\right)_{(\ell)} & = \sum_{j\not=i} c_i^{\e}c_j^{\e}  \left\{ 
\left[a_0 \frac{2\pi  m_j  \|b_{ij}\|_{L^1}}{(m_i+m_j)} +  \, a_1  \frac{10\pi kT \|b_{ij}\|_{L^1}}{m_i}\right ] \left[ (u_j)^\e_{(\ell)} - (u_i)^\e_{(\ell)} \right]
 \right. \\
+ \sum_{n\geq 2} a_n 
&  \left. \frac{ m_j (m_i m_j)^{3/2}   \|b_{ij}\|_{L^1}}{4\pi^2\e(m_i+m_j)( k\, T)^3}  \int_{\R^{6}} |v-v_*|^{2n} \left({v_*}_{(\ell)}-v_{(\ell)}\right)
e^{-{m_i}\vert v-\e u^{\e}_i \vert^2/2kT}
e^{-{m_j}\vert v_*-\e u^{\e}_j \vert^2/2kT}
\, \xd v_*\, \xd v  \right \} .
\end{align*}
{The computations that lead to the coefficients of $a_0$ can be found in \cite{bou-gre-sal-15}.
The coefficients of the other terms $a_i$ are evaluated by performing the change of variables: $(v,v_*)\mapsto (v+\e u^{\e}_i, v_*+\e u^{\e}_j)$ and then computing various Gaussian integrals (whose explicit expressions can be found in the Appendix).}

Making a change of variables: $(v,v_*)\mapsto (v+\e u^{\e}_i, v_*+\e u^{\e}_j)$ in the {integral term of the expression written above} yields:
\begin{align*}
\left(\Theta^\e_i\right)_{(\ell)} & = 
\sum_{j\not=i} c_i^{\e}c_j^{\e} \left\{ 
\left[a_0 \frac{2\pi  m_j  \|b_{ij}\|_{L^1}}{(m_i+m_j)} +  \, a_1  \frac{10\pi kT \|b_{ij}\|_{L^1}}{m_i}\right ]\left[ (u_j)^\e_{(\ell)} - (u_i)^\e_{(\ell)} \right]
\right. \\
& 
+ \sum_{n\geq 2}a_n 
 \frac{ m_j (m_i m_j)^{3/2} \|b_{ij}\|_{L^1}}{4\pi^2\e(m_i+m_j)(k\, T)^3}  
\left. \!\! \int_{\R^{6}} \!\! |v+\e u^{\e}_i-v_*-\e u^{\e}_j|^{2n} \! \left({v_*}_{(\ell)} + \e (u_j)^{\e}_{(\ell)}-v_{(\ell)} - \e (u_i^{\e})_{(\ell)}\right)
e^{-{m_i}\vert v \vert^2/2kT}
e^{-{m_j}\vert v_* \vert^2/2kT}
\, \xd v_*\, \xd v \right\}\! .
\end{align*}
{
We employ the multinomial theorem:
\[(b_1+b_2+\cdots+b_k)^n=\sum_{\substack{j_1,\ j_2,\dots,\ j_k \\ 0 \leq j_i
\leq n \text{ for each } i \\ \text{ and } j_1 + \ldots + j_k =
n}}\binom{n}{j_1, j_2, \dots , j_k}b_1^{\,\, j_1}b_2^{\,\, j_2}\cdots b_k^{\,\, j_k}\]
where the multinomial coefficients are
\begin{align*}
\binom{n}{j_1, j_2, \dots , j_k}=\dfrac{n!}{j_1! j_2!\dots j_k!}.
\end{align*}
}
This yields
\begin{align*}
 \left(\Theta^\e_i\right)_{(\ell)}  = \sum_{j\not=i} c_i^{\e}c_j^{\e} \left\{
\left[a_0 \frac{2\pi  m_j  \|b_{ij}\|_{L^1}}{(m_i+m_j)} +  \, a_1  \frac{10\pi kT \|b_{ij}\|_{L^1}}{m_i}\right ] \left[ (u_j)^\e_{(\ell)} - (u_i)^\e_{(\ell)} \right]\right.
\\
+
\sum_{n\geq 2}a_n 
 \frac{m_j (m_i m_j)^{3/2}\|b_{ij}  \|_{L^1}}{4\pi^2\e(m_i+m_j) (k\, T)^3} 
\int_{\R^{6}} \Big\{ \Big( \sum_{n_1+n_2+n_3=n}\frac{n!}{n_1! n_2! n_3!}
 \prod_{1\le r\le 3}\left({v}_{(r)} + \e (u_i^{\e})_{(r)}-v_{*(r)} - \e (u_j)^{\e}_{(r)}\right)^{2n_r}\Big)\times \\
 \left({v_*}_{(\ell)} + \e (u_j)^{\e}_{(\ell)}-v_{(\ell)} - \e (u_i^{\e})_{(\ell)}\right)
e^{-{m_i}\vert v \vert^2/2kT}
e^{-{m_j}\vert v_* \vert^2/2kT}\Big\}
\, \xd v_*\, \xd v \bigg\}.
\end{align*}
Another application of the multinomial {theorem} in the previous expression yields:
\begin{align*}
\left(\Theta^\e_i\right)_{(\ell)} = & \sum_{j\not=i}  c_i^{\e}c_j^{\e} \left\{
\left[a_0 \frac{2\pi  m_j  \|b_{ij}\|_{L^1}}{(m_i+m_j)} +  \, a_1  \frac{10\pi kT \|b_{ij}\|_{L^1}}{m_i}\right ] \left[ (u_j)^\e_{(\ell)} - (u_i)^\e_{(\ell)} \right]
\right.
\\ 
+
\sum_{n\geq 2}a_n 
 & \frac{ m_j (m_i m_j)^{3/2}
\|b_{ij}\|_{L^1}}{4\pi^2 \e(m_i+m_j)( k\, T)^3}  
\int_{\R^{6}} \Big\{ \Big( \sum_{n_1+n_2+n_3=n}\frac{n!}{n_1! n_2! n_3!} \times \\ & 
\prod_{1\le r\le 3}
\Big(\sum_{\alpha+\beta+\gamma+\lambda = 2n_r} \frac{(2n_r)!}{\alpha!\beta!\gamma!\lambda!}
\left([{v}_{(r)}]^\alpha [-v_{*(r)}]^\beta [\e (u_i^{\e})_{(r)}]^\gamma [- \e (u_j)^{\e}_{(r)}]^\lambda\right)\Big)\Big) \times \\
& \left({v_*}_{(\ell)} + \e (u_j)^{\e}_{(\ell)}-v_{(\ell)} - \e (u_i^{\e})_{(\ell)}\right)
e^{-{m_i}\vert v \vert^2/2kT}
e^{-{m_j}\vert v_* \vert^2/2kT}\Big\}
\, \xd v_*\, \xd v \bigg\}.
\end{align*}
The terms of $\mathcal{O}(\e^{-1})$ in $\left(\Theta^\e_i\right)_{(\ell)}$ are the following:
\begin{align*} 
& \sum_{j\not=i} c_i^{\e}c_j^{\e} \bigg\{
\sum_{n\geq 2}a_n
 \frac{m_j (m_i m_j)^{3/2} \|b_{ij}\|_{L^1}}{4\pi^2(m_i+m_j)( k\, T)^3} 
\int_{\R^{6}} \Big\{\Big( \sum_{n_1+n_2+n_3=n}\frac{n!}{n_1! n_2! n_3!} \times \\
& \prod_{1\le r\le 3}\left(\sum_{\alpha+\beta=2n_r} \frac{(2n_r)!}{\alpha! \beta!}\left([{v}_{(r)}]^\alpha [-v_{*(r)}]^\beta\right)\right)\Big)\left({v_*}_{(\ell)}  - v_{(\ell)} \right)
e^{-{m_i}\vert v \vert^2/2kT}
e^{-{m_j}\vert v_* \vert^2/2kT}\Big\}
\, \xd v_*\, \xd v \bigg \}.
\end{align*}
Observe that all the terms in the above sum vanish {since the integrands are odd with respect to the variables $v$ or $v_*$}. Hence there is no contribution of the terms of $\mathcal{O}(\e^{-1})$ to $\left(\Theta^\e_i\right)_{(\ell)}$. Now, we move on to consider the terms of order $\mathcal{O}(1)$ in $\left(\Theta^\e_i\right)_{(\ell)}$, which have the form:
\begin{align*}
\sum_{j\not=i} c_i^{\e}c_j^{\e} \bigg\{
 \Big( a_0 & \frac{2\pi  m_j  \|b_{ij}\|_{L^1}}{(m_i+m_j)} +  \, a_1 \frac{10\pi kT \|b_{ij}\|_{L^1}}{m_i}\Big)
 \left(  (u_j)^{\e}_{(\ell)} - (u_i^{\e})_{(\ell)} \right)
 \\
+ \sum_{n\ge2}a_n  & \frac{ m_j    (m_i m_j)^{3/2}  \|b_{ij}\|_{L^1}}{4\pi^2(m_i+m_j)(k\, T)^3}
\int_{\R^{6}} \Big\{\Big( \sum_{n_1+n_2+n_3=n}\frac{n!}{n_1! n_2! n_3!} \prod_{1\le r\le 3}\Big(\sum_{\alpha+\beta=2n_r} \frac{(2n_r)!}{\alpha! \beta!}\times \\
& \left([{v}_{(r)}]^\alpha [-v_{*(r)}]^\beta\right)\Big)\Big)\left(  (u_j)^{\e}_{(\ell)} - (u_i^{\e})_{(\ell)} \right)
e^{-{m_i}\vert v \vert^2/2kT}
e^{-{m_j}\vert v_* \vert^2/2kT}\Big\}
\, \xd v_*\, \xd v \bigg\}.
\end{align*}
Hence, we indeed have the asymptotic behaviour of $\left(\Theta^\e_i\right)_{(\ell)}$ as in \eqref{eq:lem-Theta} with the coefficients $\Delta_{ij}$ given by \eqref{eq:Delta-ij}.

\hfill$\square$ \par\medskip

\Remark {\it The kinetic collision kernel for the three dimensional hard spheres, i.e. $\Phi(|v-v_*|)=|v-v_*|$ is not an analytic function of $v-v_*$. Hence our approach cannot be directly applied to this case. However, one could approximate the hard sphere kernel by an analytic expression of the type \eqref{eq:Phi-power} and then perform the computations on the approximate series. This would yield an approximation on the binary diffusion coefficients for the hard sphere case.}
\bigskip

\subsection{Limiting behavior of the system}

Now, we are equipped to state the main result of this article. Putting together the results of Proposition \ref{prop:M1-M2} and Lemma \ref{lem:Theta}, we have indeed proved {at the formal level} the following theorem.

\begin{Thm}\label{thm:macro}
The local Maxwellian states $(\ref{ansatz})$ are solution of the initial value problem for the system
of scaled Boltzmann equations $(\ref{eqbo})$ if $(c_i^{\e},u_i^{\e})$ solves 
\begin{align}
\displaystyle
&
\pa_t c_i^{\e} +\nabla_x\cdot \left(
c_i^{\e}u_i^{\e}\right) =0, \label{euler_resc1}\\
\displaystyle
& \nabla_x
c_i^{\e} = \sum_{j\not=i} \tilde{\Delta}_{ij}\left(c^{\e}_i c^{\e}_j (u_j)^{\e}_{(\ell)} - c^{\e}_i c^{\e}_j (u_i^{\e})_{(\ell)}\right) + \mathcal{O}(\e), \label{euler_resc2}
\end{align}
{with the coefficients $\tilde{\Delta}_{ij}$ given by
\begin{equation}\label{eq:thm:Delta-ij}
\begin{array}{cl}
\displaystyle\tilde{\Delta}_{ij} = &\displaystyle a_0 \frac{2\pi  m_im_j  \|b_{ij}\|_{L^1}}{(m_i+m_j)kT} +  \, a_1 10\pi \|b_{ij}\|_{L^1}\\[0.5 cm]
& \displaystyle + \sum_{n\ge2} a_n\, \frac{2\, \pi \|b_{ij}\|_{L^1}}{k\, T} \frac{(m_i m_j)}{(m_i + m_j)}
\Big( 
\sum_{n_1+n_2+n_3=n}\frac{n!}{n_1! n_2! n_3!}\sum_{\substack{\alpha,\beta,\gamma,\delta,\rho,\eta\in 2\N^*\\ \alpha+\beta=2n_1\\ \gamma+\delta = 2n_2\\ \rho+\eta= 2n_3}} \frac{(2n_1)!}{\alpha! \beta!} \frac{(2n_2)!}{\gamma! \delta!}\frac{(2n_3)!}{\rho! \eta!}\times \\[0.3 cm]& 
\mathcal{E}(\alpha, \beta, \gamma, \delta, \rho, \eta)
\left(
\frac{kT}{m_i}
\right)^{(\alpha+\gamma+\rho)/2}
\left(
\frac{kT}{m_j}
\right)^{(\beta+\delta+\eta)/2}\Big),
\end{array}
\end{equation}
where
\begin{equation}\label{eq:mathcal-E}
\begin{aligned}
\mathcal{E}(\alpha, \beta, \gamma, \delta, \rho, \eta) := \left(
(\alpha -1)(\alpha -3)\cdots 1
\right)
\left(
(\beta -1)(\beta -3)\cdots 1
\right)
\left(
(\gamma -1)(\gamma -3)\cdots 1
\right)\times \\
\left(
(\delta -1)(\delta -3)\cdots 1
\right)
\left(
(\rho -1)(\rho -3)\cdots 1
\right)
\left(
(\eta -1)(\eta -3)\cdots 1
\right).
\end{aligned}
\end{equation}
}
\end{Thm}
\noindent{\it Proof}\par

\noindent
{
From \eqref{M2} and \eqref{eq:lem-Theta}, we have: 
\begin{align*}
\nabla_x c^\e_i
=
\sum_{j\not=i} \frac{m_i}{kT} \Delta_{ij}\left(c^{\e}_i c^{\e}_j (u_j)^{\e}_{(\ell)} - c^{\e}_i c^{\e}_j (u_i^{\e})_{(\ell)}\right) + \mathcal{O}(\e).
\end{align*}
As the coefficients $\Delta_{ij}$ involve Gaussian integrals, by using the expressions from the Appendix, we first compute the integral terms in \eqref{eq:Delta-ij}:
\begin{align*}
& \int_{\R^{6}} \Big\{\Big( \sum_{n_1+n_2+n_3=n}\frac{n!}{n_1! n_2! n_3!}\times
\prod_{1\le r\le 3}\Big(\sum_{\alpha+\beta=2n_r} \frac{(2n_r)!}{\alpha! \beta!} \left([{v}_{(r)}]^\alpha [-v_{*(r)}]^\beta\right)\Big)\Big)  \times\\
&
e^{-{m_i}\vert v \vert^2/2kT}
e^{-{m_j}\vert v_* \vert^2/2kT}\Big\}
\, \xd v_*\, \xd v\\
& = \int_{\R^{6}} \!\! \Big\{\Big( \!\! \sum_{n_1+n_2+n_3=n} \!\! \frac{n!}{n_1! n_2! n_3!} \! \sum_{\substack{\alpha,\beta,\gamma,\delta,\rho,\eta\in 2\N^*\\ \alpha+\beta=2n_1\\ \gamma+\delta = 2n_2\\ 
\rho+\eta= 2n_3}} \! \! \frac{(2n_1)!}{\alpha! \beta!} \frac{(2n_2)!}{\gamma! \delta!}\frac{(2n_3)!}{\rho! \eta!} [{v}_{(1)}]^\alpha [v_{*(1)}]^\beta [{v}_{(2)}]^\gamma [v_{*(2)}]^\delta [{v}_{(3)}]^\rho [v_{*(3)}]^\eta\Big) \times\\
& e^{-{m_i}\vert v \vert^2/2kT}
e^{-{m_j}\vert v_* \vert^2/2kT}\Big\}
\, \xd v_*\, \xd v
\end{align*}
Proceeding as in the Appendix (see \eqref{eq:compute-example}), the above expression can be computed and equals:
\begin{align*}
& \sum_{n_1+n_2+n_3=n}\frac{n!}{n_1! n_2! n_3!}\sum_{\substack{\alpha,\beta,\gamma,\delta,\rho,\eta\in 2\N^*\\ \alpha+\beta=2n_1\\ \gamma+\delta = 2n_2\\ \rho+\eta= 2n_3}} \frac{(2n_1)!}{\alpha! \beta!} \frac{(2n_2)!}{\gamma! \delta!}\frac{(2n_3)!}{\rho! \eta!}
\mathcal{E}(\alpha, \beta, \gamma, \delta, \rho, \eta)\times\\
& \left(
\frac{kT}{m_i}
\right)^{(\alpha+\gamma+\rho)/2}
\left(
\frac{kT}{m_j}
\right)^{(\beta+\delta+\eta)/2}
\left(
\frac{2\pi kT}{m_i}
\right)^{3/2}
\left(
\frac{2\pi kT}{m_j}
\right)^{3/2},
\end{align*}
where $\mathcal{E}$ is given by \eqref{eq:mathcal-E}. Substituting the above expression for the integrals in $\Delta_{ij}$ would yield an explicit expression for $\tilde{\Delta}_{ij}= m_i\Delta_{ij}/kT$:
\begin{align*}
\tilde{\Delta}_{ij} = &\displaystyle a_0 \frac{2\pi  m_im_j  \|b_{ij}\|_{L^1}}{(m_i+m_j)kT} +  \, a_1 10\pi \|b_{ij}\|_{L^1}\\[0.2 cm]
& \displaystyle + \sum_{n\ge2} a_n\, \frac{\|b_{ij}\|_{L^1}}{4\pi^2(m_i+m_j)} \frac{(m_i m_j)^{5/2}}{(kT)^4}\Big( 
\sum_{n_1+n_2+n_3=n}\frac{n!}{n_1! n_2! n_3!}\sum_{\substack{\alpha,\beta,\gamma,\delta,\rho,\eta\in 2\N^*\\ \alpha+\beta=2n_1\\ \gamma+\delta = 2n_2\\ \rho+\eta= 2n_3}} \frac{(2n_1)!}{\alpha! \beta!} \frac{(2n_2)!}{\gamma! \delta!}\frac{(2n_3)!}{\rho! \eta!}\times\\[0.3 cm]
& \mathcal{E}(\alpha, \beta, \gamma, \delta, \rho, \eta)
\left(
\frac{kT}{m_i}
\right)^{(\alpha+\gamma+\rho)/2}
\left(
\frac{kT}{m_j}
\right)^{(\beta+\delta+\eta)/2}
\left(
\frac{2\pi kT}{m_i}
\right)^{3/2}
\left(
\frac{2\pi kT}{m_j}
\right)^{3/2}
\Big)\\[0.3 cm]
= &\displaystyle a_0 \frac{2\pi  m_im_j  \|b_{ij}\|_{L^1}}{(m_i+m_j)kT} +  \, a_1 10\pi \|b_{ij}\|_{L^1}\\[0.2 cm]
& \displaystyle + \sum_{n\ge2} a_n\, \frac{2\, \pi \|b_{ij}\|_{L^1}}{k\, T} \frac{(m_i m_j)}{(m_i + m_j)}
\Big( 
\sum_{n_1+n_2+n_3=n}\frac{n!}{n_1! n_2! n_3!}\sum_{\substack{\alpha,\beta,\gamma,\delta,\rho,\eta\in 2\N^*\\ \alpha+\beta=2n_1\\ \gamma+\delta = 2n_2\\ \rho+\eta= 2n_3}} \frac{(2n_1)!}{\alpha! \beta!} \frac{(2n_2)!}{\gamma! \delta!}\frac{(2n_3)!}{\rho! \eta!}\times\\[0.3 cm]
& \mathcal{E}(\alpha, \beta, \gamma, \delta, \rho, \eta)
\left(
\frac{kT}{m_i}
\right)^{(\alpha+\gamma+\rho)/2}
\left(
\frac{kT}{m_j}
\right)^{(\beta+\delta+\eta)/2}\Big).
\end{align*}
}
\hfill$\square$ \par\medskip

Note that the coefficients $\tilde{\Delta}_{ij}$ are symmetric with respect to each pair of species since $b_{ij}=b_{ji}$ and that the structure is much more intricate than the corresponding binary diffusion coefficients computed with Maxwellian cross section in \cite{bou-gre-sal-15}.

\Remark {\it Observe that the first term in the expression for $\tilde{\Delta}_{ij}$ is nothing but the expression obtained in} \cite{bou-gre-sal-15}. {\it Note also that $\tilde{\Delta}_{ij}$ in $(\ref{eq:thm:Delta-ij})$ is given in terms of the reduced masses of the species and the temperature $T$ except for the term involving $a_1$ in the analytic expression $(\ref{eq:Phi-power})$. In particular, if the kinetic collision kernel has the form $\Phi(\vert v-v_* \vert)= \vert v-v_* \vert^2$, then the preceding computations in the paper yield the following expression for $\left(\Theta^\e_i\right)_{(\ell)}$:
\begin{align*}
\left(\Theta^\e_i\right)_{(\ell)} = \sum_{j\not=i}\frac{2\pi m_j \|b_{ij}\|_{L^1}}{(m_i+m_j)}\left(\frac{5kT}{m_i} + \frac{5kT}{m_j}\right)\left(c^{\e}_i c^{\e}_j (u_j)^{\e}_{(\ell)} - c^{\e}_i c^{\e}_j (u_i^{\e})_{(\ell)}\right) + \mathcal{O}(\e),
\end{align*}
thus giving the following expression for $\tilde{\Delta}_{ij}$:
\begin{align*}
\tilde{\Delta}_{ij} = a_1 10\pi \|b_{ij}\|_{L^1}.
\end{align*}
The rest of the terms in the expression (\ref{eq:thm:Delta-ij}), however, do depend on the temperature $T$ and the dependence is non-trivial.}
\bigskip

In the following, let us set 
$$
F_i^\e(t,x)= \frac{1}{\e}\int_{\R^3}v\, f_i^\e(t,x,v)\,\xd v=
c_i^{\e}(t,x)u_i^{\e}(t,x), \quad \mbox{ for } (t,x)\in\R_+\times\R^3, 
$$
and denote, for any $t\ge0$
and $x\in\R^3$, 
$$
c_i(t,x) = \lim_{\e\to 0^+} c_i^{\e}(t,x), \qquad F_i(t,x)= \lim_{\e\to
0^+}F_i^\e(t,x).
$$
In the limit, Equations (\ref{euler_resc1})--(\ref{euler_resc2}) give a system of
equations, which has the following form for the density-flux set of unknown
$(c_i,F_i)$:
\begin{align}
&\displaystyle
\pa_t c_i +\nabla_x\cdot F_i=0, \phantom{\int} \nonumber\\
&\displaystyle
\nabla_x c_i = -\sum_{j\neq i}\tilde{\Delta}_{ij}({c_jF_i- c_iF_j}).
\label{e:almostms2}
\end{align}
An argument similar to the ones in \cite{bou-gre-sal-15}, which consists in writing the approximation of the conservation of the kinetic energy at the leading order in $\e$ for the Boltzmann system, allows to deduce that
$$
c=\sum c_i=\sum \cini=1.
$$
We can hence obtain, in the limit, the Maxwell-Stefan system in the case of analytic cross-section and Grad's cutoff assumption:
\begin{equation*}
\left\{
\begin{array}{l}
\displaystyle
\pa_t c_i + \na_x\cdot F_i = 0 \quad \mbox{on } \R_+\times \R^3,\\ \\
\displaystyle -c\, \na_x n_i = \frac 1c\sum_{j\neq i} \frac{c_j F_i-c_i F_j}{\dij} \quad
\mbox{on } \R_+^*\times \R^3,
\end{array}
\right.
\end{equation*}
where the binary diffusion coefficient have the form 
{
\begin{align*}
\dij
= \frac 1 c \displaystyle \Bigg \{a_0 \frac{2\pi  m_im_j  \|b_{ij}\|_{L^1}}{(m_i+m_j)kT} +  \, a_1 10\pi \|b_{ij}\|_{L^1}
+ \sum_{n\ge2} a_n\, \frac{2\, \pi \|b_{ij}\|_{L^1}}{k\, T} \frac{(m_i m_j)}{(m_i + m_j)}
\Big( 
\sum_{n_1+n_2+n_3=n}\frac{n!}{n_1! n_2! n_3!}
 \times
\\
\sum_{\substack{\alpha,\beta,\gamma,\delta,\rho,\eta\in 2\N^*\\ \alpha+\beta=2n_1\\ \gamma+\delta = 2n_2\\ \rho+\eta= 2n_3}} \frac{(2n_1)!}{\alpha! \beta!} \frac{(2n_2)!}{\gamma! \delta!}\frac{(2n_3)!}{\rho! \eta!}\mathcal{E}(\alpha, \beta, \gamma, \delta, \rho, \eta)
\left(
\frac{kT}{m_i}
\right)^{(\alpha+\gamma+\rho)/2}
\left(
\frac{kT}{m_j}
\right)^{(\beta+\delta+\eta)/2}\Big)\Bigg \}^{-1}.
\end{align*}
}

\section{Appendix}\label{sec:explicit}
{
The objective of this appendix is the explicit computation of the Gaussian integrals appearing in the expression for ${\Delta}_{ij}$ given by \eqref{eq:Delta-ij}
leading to the expression of $\tilde{\Delta}_{ij}$ given by \eqref{eq:thm:Delta-ij}. We have the following normalization:
\begin{align*}
\int_{\R^3}
\left(
\frac{m_i}{2\pi k T}
\right)^{3/2}
e^{-m_i|v|^2/2kT}
\, {\rm d}v
= 1.
\end{align*}
We observe that
\begin{align*}
& \left(
\frac{m_i}{2\pi k T}
\right)^{3/2}
\int_{\R^3}
[v_{(1)}]^\alpha\, 
[v_{(2)}]^\beta\, 
[v_{(3)}]^\gamma\, 
e^{-m_i|v|^2/2kT}
\, {\rm d}v_{(1)}\, {\rm d}v_{(2)}\, {\rm d}v_{(3)}\\
& = -\left(
\frac{m_i}{2\pi k T}
\right)^{3/2}
\int_{\R^3}
\left(
\frac{kT}{m_i}
\right)
[v_{(1)}]^{\alpha - 1}\, 
\partial_{v_{(1)}} e^{-m_i|v|^2/2kT}\, 
[v_{(2)}]^\beta\, 
[v_{(3)}]^\gamma\, 
\, {\rm d}v_{(1)}\, {\rm d}v_{(2)}\, {\rm d}v_{(3)}.
\end{align*}
Performing integration by parts with respect to $v_{(1)}$ variable we arrive at
\begin{align*}
& \left(
\frac{m_i}{2\pi k T}
\right)^{3/2}
\int_{\R^3}
\left(
\frac{kT}{m_i}
\right)
(\alpha - 1)
[v_{(1)}]^{\alpha - 2}\, 
e^{-m_i|v|^2/2kT}\, 
[v_{(2)}]^\beta\, 
[v_{(3)}]^\gamma\, 
\, {\rm d}v_{(1)}\, {\rm d}v_{(2)}\, {\rm d}v_{(3)}\\
& = -\left(
\frac{m_i}{2\pi k T}
\right)^{3/2}
\int_{\R^3}
\left(
\frac{kT}{m_i}
\right)^2
(\alpha - 1)
[v_{(1)}]^{\alpha - 3}\, 
\partial_{v_{(1)}} e^{-m_i|v|^2/2kT}\, 
[v_{(2)}]^\beta\, 
[v_{(3)}]^\gamma\, 
\, {\rm d}v_{(1)}\, {\rm d}v_{(2)}\, {\rm d}v_{(3)}.
\end{align*}
Performing an integration by parts in the $v_{(1)}$ variable again yields:
\begin{align*}
\left(
\frac{m_i}{2\pi k T}
\right)^{3/2}
\int_{\R^3}
\left(
\frac{kT}{m_i}
\right)^2
(\alpha - 1)
(\alpha - 3)
[v_{(1)}]^{\alpha - 4}\, 
e^{-m_i|v|^2/2kT}\, 
[v_{(2)}]^\beta\, 
[v_{(3)}]^\gamma\, 
\, {\rm d}v_{(1)}\, {\rm d}v_{(2)}\, {\rm d}v_{(3)}
\end{align*}
Continuing in a similar fashion with respect to all the variables, i.e. $v_{(1)}, v_{(2)}, v_{(3)}$, we arrive at the following expression:
\begin{equation}\label{eq:compute-example}
\begin{aligned}
& \left(
\frac{m_i}{2\pi k T}
\right)^{3/2}
\int_{\R^3}
[v_{(1)}]^\alpha\, 
[v_{(2)}]^\beta\, 
[v_{(3)}]^\gamma\, 
e^{-m_i|v|^2/2kT}
\, {\rm d}v_{(1)}\, {\rm d}v_{(2)}\, {\rm d}v_{(3)}\\
& = \left(
(\alpha -1)(\alpha -3)\cdots 1
\right)
\left(
(\beta -1)(\beta -3)\cdots 1
\right)
\left(
(\gamma -1)(\gamma -3)\cdots 1
\right)
\left(
\frac{kT}{m_i}
\right)^{(\alpha+\beta+\gamma)/2}.
\end{aligned}
\end{equation}
}

\medskip
\medskip

\noindent {\bf Acknowledgments:} The authors thank the referees for their useful comments and remarks.
The authors appreciate the fruitful suggestions of Claude Bardos during the preparation of this article.
This work was partially funded by the projects \textit{Kimega} (ANR-14-ACHN-0030-01) and
\textit{Kibord} (ANR-13-BS01-0004).
H.H. acknowledges the support of the ERC grant \textit{MATKIT}.
F.S. would like to thank the Isaac Newton Institute for Mathematical Sciences, Cambridge, for support and hospitality during the programme \textit{Mathematical Modelling and Analysis of Complex Fluids and Active Media in Evolving Domains}. This work was partially supported by EPSRC grant no EP/K032208/1.

\end{document}